\documentclass[11pt]{article}
\usepackage{graphicx, colordvi}
\usepackage{amsfonts}
\usepackage{amssymb}
\usepackage{amsthm,cite,color}
\usepackage{dsfont}
\usepackage{epsfig}
\usepackage{mathrsfs}
\usepackage{amsfonts}
\usepackage{amssymb}
\usepackage{amssymb,amsfonts,amsmath,amsthm,cite,color}
\usepackage{dsfont}
\usepackage{epsfig}
\usepackage{mathrsfs}
\usepackage{longtable}

\newtheorem{theo}{Theorem}
\newtheorem{prop}[theo]{Proposition}

\newtheorem{defi}[theo]{Definition}
\newtheorem{lem}[theo]{Lemma}
\newtheorem{exm}[theo]{Example}

\makeatletter \@addtoreset{equation}{section}
\@addtoreset{theo}{section} \makeatother

\newcommand{\bQ} { {\mathbb{Q}}}
\newcommand{\bZ} { {\mathbb{Z}}}

\def\qed{\hfill \rule{4pt}{7pt}}
\def\pf{\noindent {\it Proof.} }

\begin{document}
\begin{center}
{\large\bf Infinite Orders and Non-$D$-finite Property \\[5pt]
 of $3$-Dimensional Lattice Walks}

\vskip 3mm

Daniel K. Du$^1$, Qing-Hu Hou$^2$ and Rong-Hua Wang$^3$\\[7pt]

$^{1,2}$Center for Applied Mathematics \\
Tianjin University, Tianjin 300072, P. R. China \\[5pt]

$^{3}$Center for Combinatorics, LPMC-TJKLC \\
Nankai University, Tianjin 300071, P. R. China

\vskip 3mm

 E-mail: $^1$daniel@tju.edu.cn, $^2$hou@nankai.edu.cn,
 $^3$wangwang@mail.nankai.edu.cn
\end{center}

\begin{abstract}
Recently, Bostan and his coauthors investigated lattice walks restricted to the non-negative octant $\mathbb{N}^3$.
For the $35548$ non-trivial models with at most six steps, they found that many models associated to a group
of order at least $200$ and conjectured these groups were in fact infinite groups.
In this paper, we first confirm these conjectures and then consider the non-$D$-finite property of the generating function for some of these models.
\end{abstract}

\section{Introduction}
The objective of this paper is to use the properties of Jacobian matrix at fixed points to derive the infiniteness of groups associated with certain  lattice walks restricted to the positive octant.
Furthermore, we present the non-D-finiteness of corresponding generating functions for some lattice walks of infinite order by considering the asymptotic behavior of theirs coefficients.

Counting walks in a fixed region of the lattice $\bZ^d$ is a classical topic
in enumerative combinatorics~\cite{Maher1978,GesselZeilberger1992,Mishna2009,Raschel2012} and in probability theory~\cite{Narayana1979,Mohanty1979}.
In the past few years, lattice path models restricted to the quarter plane and the positive octant have received special attention, and recent works~\cite{BBKM2015,MelouMishna2010,BostanKauers2010,KurkovaRaschel2011
,BostanRaschelSalvy2014,MishnaRechnitzer2009,MelczerMishna2014} have shown how they can help us better understand generating functions of lattice walks.

Many recent papers dealt with the enumeration of lattice walks with prescribed steps confined to the positive quadrant. In fact, Bousquet-m\'{e}lou and Mishna~\cite{MelouMishna2010} proved that among the $2^8$ possible cases of
small-step in the quarter plane, there were exactly $79$ inherently different cases.
Then, they showed that $23$ of these models were associated with finite group, of which $22$ ones admitted $D$-finite generating functions (see, for example~\cite{Stanley1980} for an overview on D-finite) .  The $23$rd model, known as Gessel walks, was proven $D$-finite, and even algebraic, by Bostan and Kauers~\cite{BostanKauers2010}.
Moreover, it was conjectured in \cite{MelouMishna2010} that the $56$ remaining models with infinite group had non-$D$-finite generating functions. This was proved by  Kurkova and Raschel~\cite{KurkovaRaschel2011} for the $51$ nonsingular walks. The remaining $5$ singular models were proven by Mishna and Rechnitzer~\cite{MishnaRechnitzer2009} and Melczer and Mishna~\cite{MelczerMishna2014}. The classification is now complete for walks with steps in $\{0,\pm1\}^2$: the generating function is D-finite if and only if a certain group associated with the model is finite.

Recently, Bostan and his coauthors~\cite{BBKM2015} considered the analogous problem for lattice walks confined to the non-negative octant $\mathbb{N}^3$. They
showed there were $35548$ non-trivial models with at most six steps. Each model corresponds to a group which plays an important role in exploring the properties of
the generating function. They found that many models associated to a group
of order at least $200$ and conjectured these groups were in fact infinite groups.

In this paper, we mainly utilize two methods employed by  Bousquet-M\'{e}lou and Mishna in~\cite{MelouMishna2010} to confirm these conjectures by considering models of dimension two and three, respectively. For the notation of dimension of a model, one can refer to  Definition \ref{DEFI:dimension}.

More specifically, for the cases of models of dimension two, Bostan $et\ al.$ \cite{BBKM2015} showed that there were $527$ models of cardinality at most $6$.
They found that $118$ models associated to a finite group of order at most $8$, and conjectured that the remaining $409$ ones associated to a group of infinite order. Our first result is to confirm this conjecture as follows.

\begin{theo}\label{CJ:409}
The $409$ two-dimensional models associated to groups of order at least $200$ are in fact associated to infinite groups.
\end{theo}

Indeed, most of these models have the property of non-$D$-finite, which means that their generating functions
 do not satisfy any non-trivial linear differential recurrences with polynomial coefficients.
\begin{theo}\label{ND:384}
For these $409$ two-dimensional models associated to infinite groups,
the generating functions of the excursions of the $366$ non-singular models are all non-$D$-finite,
and there are $18$ singular models with non-$D$-finite generating functions.
\end{theo}

For the cases of three-dimensional models, Bostan $et\ al.$ showed that there were $20634$ models associated with a group of order at least $200$  and conjectured the order to be infinite in~\cite{BBKM2015}. Our third result is to confirm this conjecture.

\begin{theo}\label{CJ:20634}
The $20634$ three-dimensional models associated with groups of order at least $200$ are in fact associated with infinite groups.
\end{theo}

This paper is organized as follows. We first recall some notations in Section $2$.
Then we derive the infiniteness of groups associated with certain models in Section $3$. Meanwhile
the proof of Theorem \ref{CJ:409} and Theorem \ref{CJ:20634} will be presented, respectively.
Section $4$ discusses the non-$D$-finite property and the proof of Theorem \ref{ND:384} will be presented.

\section{Preliminaries}
To make this paper self-contained, we now recall some definitions and notations.
In particular, we shall use the dimension, the characteristic polynomial and the associated group of models.

Given the hyper cubic lattice $\mathbb{Z}^d$, a finite set of \emph{steps}
$\mathcal{S}\subset\mathbb{Z}^d $ is called a \emph{model} as adopted in \cite{BBKM2015}.
We define an $\mathcal{S}$-\emph{walk} to be any walk which
starts from the origin $(0,0,0)$ and takes its steps in $\mathcal{S}$.
In particular, we focus on \emph{octant walks},
which are $\mathcal{S}$-walks remaining in the non-negative octant $\mathbb{N}^3$, with $\mathbb{N}=\{0,1,2,\dots\}$.
Then we have
\begin{defi}\label{DEFI:genfun}
The $\emph{complete generating function}$ of an octant walk is
\[
O(x,y,z;t)=\sum_{i,j,k,n\geq 0}o(i,j,k;n)x^iy^jz^kt^n,
\]
where $o(i,j,k;n)$ is the number of $n$-step walks in the octant that end at position $(i,j,k)$.
The specialization $O(0,0,0;t)$ counts $\mathcal{S}$-walks returning to the origin, called $\mathcal{S}$-\emph{excursions}.
\end{defi}
To shorten notation,
we denote steps of $\mathbb{Z}^d$ by
$d$-letter words. For example, $\overline{1}10$ stands for the step $(-1,1,0)$.
In fact, an $\mathcal{S}$-\emph{walk} of length $n$ can be viewed as a word $w=w_1w_2\cdots w_n$
made up of letters of $\mathcal{S}$. For each step $s\in \mathcal{S}$, let $a_s$ be the number of occurrences of $s$ in $w$.
Then $w$ ends in the positive octant if and only if the following three linear inequalities hold
\begin{equation}\label{EQ:dimension}
\sum_{s\in \mathcal{S}} a_ss_x\geq0,\ \sum_{s\in \mathcal{S}} a_ss_y\geq0,\ \sum_{s\in \mathcal{S}} a_ss_z\geq0,
\end{equation}
where $s=\{s_x,s_y,s_z\}$ are steps in $\mathcal{S}$. Furthermore, $w$ is an octant walk if the multiplicities observed in each of its prefixes satisfy these inequalities.
More generally, we give the definition of dimension of a model as follows.

\begin{defi}\label{DEFI:dimension}
Let $d\in\{0,1,2,3\}$. A model $\mathcal{S}$ is said to have dimension at most $d$ if there exist $d$ inequalities in Equation~\eqref{EQ:dimension} such that any $|\mathcal{S}|$-tuple $(a_s)_{s\in \mathcal{S}}$ of non-negative integers satisfying these $d$ inequalities satisfies in fact the three ones.
\end{defi}

Given a model $\mathcal{S}$ of cubic lattice, we denote by $S(x,y,z)$ the Laurent polynomial
\[
S(x,y,z)=\sum_{ijk\in S}x^iy^jz^k.
\]
According to the degrees of $x,y$ and $z$, respectively, $S(x,y,z)$ can be written as
\begin{align*}
S(x,y,z)& =\overline{x}A_{-}(y,z)+A_{0}(y,z)+xA_{+}(y,z)\\
        & =\overline{y}B_{-}(x,z)+B_{0}(x,z)+yB_{+}(x,z)\\
        & =\overline{z}C_{-}(x,y)+C_{0}(x,y)+zC_{+}(x,y),
\end{align*}
where $\overline{x}=1/x,\ \overline{y}=1/y,\hbox{ and }\overline{z}=1/z$. We call $S(x,y,z)$ the \emph{characteristic polynomial} of $\mathcal{S}$.

Let first assume that $\mathcal{S}$ is of $3$-dimensional. Then it has a positive step in each direction, and $A_{+}, B_{+}$ and $C_{+}$ are non-zero. Now we introduce the notation of groups associated with $\mathcal{S}$ as follows.
\begin{defi}\label{DEFI:groups}
For a given model $\mathcal{S}$, the \emph{group} associated with $\mathcal{S}$ is defined as the group $G(\mathcal{S})$ of birational transformations of the variables $[x,y,z]$ generated by the following three involutions
\[
\phi([x,y,z])=\left[\overline{x}\frac{A_{-}(y,z)}{A_{+}(y,z)},y,z\right],
\]
\[
\psi([x,y,z])=\left[x,\overline{y}\frac{B_{-}(x,z)}{B_{+}(x,z)},z\right],
\]
\[
\tau([x,y,z])=\left[x,y,\overline{z}\frac{C_{-}(x,y)}{C_{+}(x,y)}\right].
\]
\end{defi}
By construction, $G(\mathcal{S})$ fixes the characteristic polynomial $S(x,y,z)$.

For a $2$-dimensional model $\mathcal{S}$, the $z$-condition can be ignored, and the corresponding group
 $G({\cal S})$ is the group generated by $\phi$ and $\psi$.

\section{Infiniteness of Associated Groups }\label{SE:NonInfite}

In this section, we consider the $35548$ non-trivial models with at most six steps confined to the non-negative octant $\mathbb{N}^3$.
We derive the infiniteness of these groups by giving the proofs of Theorem \ref{CJ:409} and Theorem \ref{CJ:20634}, respectively.

\subsection{The Proof of Theorem \ref{CJ:409}}\label{SE:2D}
When dealing with models of dimensional two, we consider the projection of the model to a plane throughout this paper. Then the models are a multi-set of $\{\overline{1},0,1\}^2\setminus \{0,0\}$.

In order to show the infiniteness of groups associated to two dimensional octant models, we first introduce the method of fixed point argument given by Bousquet-M\'{e}lou and Mishna~\cite{MelouMishna2010} and give some preliminaries.

Assume that $\theta=\psi\circ\phi$ is well-defined in the  neighborhood of $(a,b)\in\mathbb{C}^2$, which was fixed by $\theta$.  Note that $a$ and $b$ are algebraic over $\bQ$. Let us write $\theta=(\theta_1,\theta_2)$, where $\theta_1$ and $\theta_2$ are the two coordinates of $\theta$. Each $\theta_i$ sends the pair $(x,y)$ to a rational function of $x$ and $y$. The local expansion of $\theta$ around $(a,b)$ reads
\[
\theta(a+u,b+v)=(a,b)+(u,v)J_{\theta}+O(u^2)+O(v^2)+O(uv),
\]
where $J_{\theta}$ is the Jacobian matrix of $\theta$ at $(a,b)$:
\[
J_{\theta}=\left(
    \begin{matrix}{}
       \frac{\partial \theta_1}{\partial x}(a,b) & \frac{\partial \theta_2}{\partial x}(a,b) \\[6pt]
       \frac{\partial \theta_1}{\partial y}(a,b) & \frac{\partial \theta_2}{\partial y}(a,b)
    \end{matrix}
    \right).
\]
Iterating the above expansion gives, for $m\geq 1$,
\begin{equation}\label{EQ:J^m}
\theta^m(a+u,b+v)=(a,b)+(u,v)J_{\theta}^m+O(u^2)+O(v^2)+O(uv).
\end{equation}
Assume that $\theta$ is of order $n$. Then $\theta^n(a+u,b+v)=(a,b)+(u,v)$ and Equation~\eqref{EQ:J^m} show that $J_{\theta}^n$ is the identity matrix. In particular, all eigenvalues of $J_{\theta}$ are roots of unity. This provides a strategy for proving that a group $G(\mathcal{S})$ is infinite.

We now give some properties on the fixed points of $\theta=\psi \circ \phi$ and the Jacobian matrices, which will simplify our computations.
\begin{prop}
$(a,b)$ is a fixed point of $\theta$ if and only if it is a fixed point of $\phi$ and $\psi$.
\end{prop}
\pf
Suppose $(a,b)$ is a fixed point of $\theta$. Assume that $\phi(a,b)=(u,b)$. Then we have $\psi(u,b)=(a,b)$. By definition, $\psi$ preserves the first coordinate. We thus have $u=a$ and $(a,b)$ is a fixed point of $\phi$ and $\psi$. The inverse assertion holds straightforwardly. \qed

This proposition indicates that the fixed point $(a,b)$ of $\theta$ can be determined by the equations
\[
  \frac{A_{-}(b)}{A_{+}(b)} = a^2  \quad \mbox{and} \quad \frac{B_{-}(a)}{B_{+}(a)} = b^2.
\]
Moreover, we require that $a$ and $b$ are both non-zero. Now we rewrite the left hand sides of the above two equations in reduced form by canceling the common divisor of the numerator and the denominator and we get
\[
  \frac{p_1(b)}{q_1(b)} = a^2  \quad \mbox{and} \quad \frac{p_2(a)}{q_2(a)} = b^2.
\]
We need to find the solutions of the polynomial systems
\[
  p_1(y) - x^2 q_1(y) = 0, \quad p_2(x) - y^2 q_2(x) = 0, \quad xy \not = 0.
\]
The command {\tt RegSer} in Maple package {\tt epsilon} by D.~M. Wang~\cite{Wang2002} can solve such system. By using
\begin{center}
{\tt
  RegSer([[$p_1(y) - x^2 q_1(y), p_2(x) - y^2 q_2(x)$], [$xy$]], [$x,y$]),
}
\end{center}
we will obtain a basis on the equations satisfied by the fixed points.
When the output is {\tt []}, there is no fixed points and the method fails.


The determinant of the Jacobian matrix $J_\theta$ at fixed points satisfies the following property.
\begin{lem}
The determinant of the Jacobian matrix $J_\theta$ at fixed points is~$1$.
\end{lem}
\pf By the chain rule, we have
$J_\theta = J_\psi \cdot J_\phi$. While
\[
J_\phi = \left( \begin{array}{cc} -\frac{1}{a^2} \frac{p_1(b)}{q_1(b)} & \left. \frac{\partial (p_1(y)/xq_1(y))}{\partial y} \right|_{(a,b)} \\ 0 & 1  \end{array} \right)
 = \left( \begin{array}{cc} -1 & * \\ 0 & 1  \end{array} \right),
\]
and
\[
J_\psi = \left( \begin{array}{cc} 1 & 0 \\
\left. \frac{\partial (p_2(x)/yq_2(x))}{\partial x} \right|_{(a,b)} & -\frac{1}{b^2} \frac{p_2(a)}{q_2(a)} \end{array} \right)
 = \left( \begin{array}{cc} 1 & 0 \\ * & -1  \end{array} \right).
\]
Therefore, the determinant of $J_\theta$ is $(-1) \cdot (-1)=1$. \qed

Let $p(X,x,y)$ be the numerator of
\[
\chi(X)=det(X Id-J_{\theta})=X^2 - \left( \frac{\partial (p_2(x)/yq_2(x))}{\partial x} \cdot \frac{\partial (p_1(y)/xq_1(y))}{\partial y} - 2 \right) X + 1.
\]
Once again, we use
\begin{center}
  \tt RegSer([[$p(X,x,y), p_1(y) - x^2 q_1(y), p_2(x) - y^2 q_2(x)$], [$xy$]], [$X,x,y$])
\end{center}
to obtain an equation $q(X)$ satisfied by $X$, the eigenvalues of $J_\theta$. To verify whether the eigenvalues of $J_\theta$ are roots of unit, we need only to check whether all the irreducible factors of $q(X)$ are cyclotomic polynomials.

To make the above statements easier understood, we present two examples.
\begin{exm}
Suppose $\mathcal{S}=[\overline{1}0, \overline{1}1, \overline{1}1, \overline{1}1,0\overline{1}, 1]$, then the corresponding characteristic polynomial is
\[
S(x,y)=\frac{1}{x}+\frac{3y}{x}+\frac{1}{y}+xy,
\]
and
\[
A_{-}(y)=1+3y,\ A_{+}(y)=y,\ B_{-}(x)=1,\ B_{+}(x)=3/x+x.
\]
Applying the command {\tt RegSer}, we find that the fixed point of $\theta$ must satisfy the following two equations:
\[
9x-6x^3+x^5-3-x^2=0 \text{ and } -1-3y+x^2y=0.
\]
Let $p(X,x,y)$ be the numerator of $\chi(X)$, we get
\begin{align*}
p(X,x,y)=&9X^2y^3x+6X^2y^3x^3+X^2y^3x^5+3X-Xx^2+18Xxy^3\\
&+12Xx^3y^3+2Xx^5y^3+9xy^3+6x^3y^3+x^5y^3.
\end{align*}
Using {\tt RegSer} once again, we find $X$ satisfies $q(X)=0$, where
\begin{align*}
q(X)=& 27X^{10}-216X^9-2267X^8-7881X^7-15249X^6-18785X^5\\
     & -15249X^4-7881X^3-2267X^2-216X+27.
\end{align*}
It's easy to check that $q(X)$ has two irreducible factors
$
X^2+X+1
$
and
~$
27X^8-243X^7-2051X^6-5587X^5-7611X^4-5587X^3-2051X^2-243X+27.
$
Since the second factor is not a cyclotomic polynomial, we conclude that~$\mathcal{S}=[\overline{1}0, \overline{1}1, \overline{1}1, \overline{1}1,0\overline{1}, 1]$ is associated with an infinite group.
\end{exm}

\begin{exm}\label{EX:not determined}
Suppose $\mathcal{S}=[\overline{1}1, \overline{1}1, 1\overline{1}, 10]$, then the corresponding characteristic polynomial is
\[
S(x,y)=\frac{2y}{x}+\frac{x}{y}+x,
\]
and
\[
A_{-}(y)=2y,\ A_{+}(y)=1/y+1,\ B_{-}(x)=x,\ B_{+}(x)=2/x.
\]
Applying the command {\tt RegSer}, the output is {\tt []} and the method fails.
\end{exm}
By this method, we show that neither eigenvalues of $379$ models are roots of unity and hence $\theta$ is an element of infinite order. There are $30$ models left, such as $\mathcal{S}=[\overline{1}1, \overline{1}1, 1\overline{1}, 10]$ in Example~\ref{EX:not determined}. Canceling the repeated steps, all these models fall in the five models (or their $x/y$ reflection) which had been proved associated with an infinite group by the valuation argument.

The valuation argument was given by Bousquet-M\'{e}lou and Mishna in~\cite{MelouMishna2010}. In fact, they defined the \emph{valuation} of a Laurent series $F(t)$ to be the smallest $d$ such that $t^d$ occurs in $F(t)$ with a non-zero coefficient. Suppose $z$ is an indeterminate, and $x$, $y$ are Laurent series in $z$ with coefficients in $\bQ$, of respective valuations $a$ and $b$. Assuming that the \emph{trailing} coefficients of these series, namely $[z^a]x$ and $[z^b]y$, are positive. Defining $x'$ by $\phi(x,y)=(x',y)$. Then the trailing coefficient of $x'$ (and $y$) is positive, and it's easy to check that the valuation of $x'$ (and $y$) only depends on $a$ and $b$:
\[
 \Phi(a,b):=(\hbox{val}(x'),\hbox{val}(y))=\begin{cases}
  {\left(-a+b(v_{-1}^{(y)}-v_{1}^{(y)}),b\right)}, & \mbox{if $b \ge 0$}; \\[15pt]
   \left(-a+b(d_{-1}^{(y)}-d_{1}^{(y)}),b\right), & \mbox{if $b \leq 0$};
\end{cases}
\]
where $v_i^{(y)}$ (resp. $d_{i}^{(y)}$) denotes the valuation (resp. degree) in $y$ of $A_i(y)$, for $i=\pm1$.
Similarly, $\psi(x,y)=(x,y')$ is well defined, and the valuations of $x$ and $y'$ only depend on $a$ and $b$:
\[
 \Psi(a,b):=(\hbox{val}(x),\hbox{val}(y'))=\begin{cases}
  {\left(a,-b+a(v_{-1}^{(x)}-v_{1}^{(x)})\right)}, & \mbox{if $a \ge 0$}; \\[15pt]
   \left(a,-b+a(d_{-1}^{(x)}-d_{1}^{(x)})\right), & \mbox{if $a \leq 0$};
\end{cases}
\]
where $v_i^{(x)}$ (resp. $d_{i}^{(x)}$) denotes the valuation (resp. degree) in $x$ of $B_i(x)$, for $i=\pm1$. For a given model $\mathcal{S}$, in order to prove the associated group $G$ is infinite, it suffices to prove that the group $G'$ generated by $\Phi$ and $\Psi$ is infinite. To prove the latter statement, it suffices to exhibit $(a,b)\in\mathbb{Z}^2$ such that the orbit of $(a,b)$ under the action of $G'$ is infinite. For the five singular models, Bousquet-M\'{e}lou and Mishna derived by induction on $n$ that
\[
(\Psi\circ\Phi)^n(1,2)=(2n+1,2n+2)\quad\hbox{ and }\quad
\Phi(\Psi\circ\Phi)^n=(2n+3,2n+2).
\]
Hence the orbit of $(1,2)$ under the action of $\Phi$ and $\Psi$ is infinite, and so are the groups $G'$ and $G$.

It's easy to check that the repeated steps do not change the value of $\Phi(a,b)$ and $\Psi(a,b)$ by the definition.
Thus, we obtain the fact that and the left $30$ models are associated with infinite groups.

This completes the proof of Theorem \ref{CJ:409}.

\subsection{The Proof of Theorem \ref{CJ:20634}}\label{SE:3D}
In this section, we consider the three-dimensional models. The proof of Theorem \ref{CJ:20634} is similar to the proof for the two dimension case.
Indeed, for three dimension cases, we could consider $\phi \circ \psi$, $\phi \circ \tau$ and $\psi \circ \tau$,
instead of $\theta = \psi \circ \phi$ in the cases of two-dimensional models. Moreover, we need only to concern two variables by fixing the third variable with any given value. For simplicity, we set the third variable to be $1/7$.

The following lemma indicates that we need only to consider one of $\phi \circ \psi$ and $\psi \circ \phi$.
\begin{lem}
If the eigenvalues of $J_{\phi \circ \psi}$ are roots of unit, then so are $J_{\psi \circ \phi}$.
\end{lem}
\pf Notice that the determinant of $J_\theta$ is $1$. The eigenvalues of $J_\theta$ are both roots of unit or neither of the eigenvalues is root of unit. Since $(\psi \circ \phi)^{-1} = \phi \circ \psi$, we have
$J_{\phi \circ \psi}^{-1} =  J_{\psi \circ \phi}$. Thus, the eigenvalues of $J_{\phi \circ \psi}$ are the reciprocal of those of $J_{\psi \circ \phi}$ and hence they are both roots of unit or none of them are roots of unit. \qed

By the fixed point method just as in Section \ref{SE:2D}, we are left $69$ models that can not be proved to be infinity. By projecting these models to two dimension models (we have three choices) and remove the repeated steps, one can find that they all fall in the five models which have been proved with an infinity group by the valuation argument. Thus the left $69$ are all associated with infinite groups.

This completes the proof of Theorem \ref{CJ:20634}.

\section{The non-D-finite Property}

In this section, we mainly discuss the non-$D$-finite property of the generating function of
the $409$ two-dimensional models associated with an infinite order, by giving the proof of
Theorem \ref{ND:384}.

As shown in Section~\ref{SE:2D}, by projected to a plane, these two-dimensional models are reduced to multi-sets of $\{\overline{1},0,1\}^2\setminus \{0,0\}$.
For a $2$D octant model where the $z$-condition is redundant, we focus on the complete generating function
\begin{equation}\label{EQ:genfun2D}
O(x,y;t):=O(x,y,1;t),
\end{equation}
which counts quadrant walks with steps in multiset $S'=\{ij:ijk\in S\}$.
The main objective of this section is to study
the non-$D$-finite property of $O(x,y;t)$.

Firstly, we consider the $\emph{nonsingular walks}$, that is for walks having at least one step from the set $\{(-1,0), (-1,-1), (0,-1)\}$. Bostan $ et\ al.$ proved that the excursion corresponding to any of the $51$ nonsingular models having no repeated step and with infinite group were not D-finite in~\cite{BostanRaschelSalvy2014}. They utilized the fact that, in many cases, we can detect non-D-finiteness of power series by looking at the asymptotic behavior of its coefficients, which is a consequence of the theory of $G$-functions and provided the following theorem.
\begin{theo}\label{TH:NotD}
Let $(a_n)_{n\geq 0}$ be an integer-valued sequence whose $n$-th term $a_n$ behaves
asymptotically like $K\cdot \rho^n\cdot n^{\alpha}$, for some real constant $K>0$. If the growth constant $\rho$ is transcendental, or if the singular exponent $\rho$ is irrational, then the generating function $A(t)=\sum_{n>0}a_nt^n$ is not D-finite.
\end{theo}

Bostan $ et\ al.$ considered the non-degeneracy of the walk: for all $(i,j)\in \mathbb{N}^2$, the set $\{n\in\mathbb{N} \colon o(i,j;n)\neq 0\}$ is nonempty; furthermore, the walk is said to be $\emph{aperiodic}$ when the gcd of the elements of this set is $1$ for all $(i,j)$. Otherwise, it is $\emph{periodic}$ and this gcd is the period. Then they restated a result of Denisov and Wachtel~\cite{DenisovWachtel2013} in the following way that can be used directly in our computations.

\begin{theo}\label{TH:asymptotic}
Let $\mathcal{S}\subset\{0,\pm1\}^2$ be the step set of a walk in the quarter plane $\mathbb{N}^2$, which is not contained in a half-plane. Let $e_n$ denote the number of excursions of length $n$ using only steps in $\mathcal{S}$, and let $\chi$ denote the characteristic polynomial $\sum_{(i,j)\in\mathcal{S}}x^iy^j$ of the step set $\mathcal{S}$. Then the system
\[
\frac{\partial\chi}{\partial x}=\frac{\partial\chi}{\partial y}=0
\]
has a unique solution $(x_0,y_0)\in\mathbb{R}_{>0}^2$. Next, define
\[
\rho:=\chi(x_0,y_0),\quad c:=\frac{\frac{\partial^2\chi}{\partial x\partial y}}{\sqrt{\frac{\partial^2\chi}{\partial^2x}\cdot\frac{\partial^2\chi}{\partial^2\partial y^2}}}(x_0,y_0),\quad \alpha:=-1-\frac{\pi}{\arccos(-c)}.
\]
Then there exists a constant $K>0$, which only depends on $\mathcal{S}$, such that
\begin{itemize}
\item
if the walk is aperiodic, then
$e_n\sim K\cdot\rho^n\cdot n^{\alpha}.$
\item
if the walk is periodic(then of periodic $2$), then
\[e_{2n}\sim K\cdot\rho^{2n}\cdot(2n)^\alpha,\quad e_{2n+1}=0.\]
\end{itemize}
\end{theo}

Then they gave an algorithmic proof that for any of the $51$ nonsingular models confined to the positive quadrant, the singular exponent $\alpha$ in the asymptotic expansion of excursion sequence was an irratinal number. Thus by the above two Theorems, the generating function $O(0,0;t)$ is not D-finite.

We note that Theorem~\ref{TH:asymptotic} still holds for multi-sets, since the repetition of a step just change the probability of the appearance of this step. Then we can apply the algorithmic irrational proof, given in Section~$2.4$ in~\cite{BostanRaschelSalvy2014}, to the $409$ two dimensional models associated to groups of infinite order. It turns out that the singular exponent $\alpha$ is irrational for $366$ nonsingular models, which proves that the corresponding excursion generating function $O(0,0;t)$ is not D-finite for these models.

The algorithmic irrational proof fails for the $43$ singular models, which were listed in the Appendix, Table~\ref{TB:43}. We find that all these models can be reduced to one of the $5$ singular step sets or their $x/y$ symmetry in two dimensional walks, when get rid of repeated steps. The $5$ singular models were proven with non-D-finite generating function by Mishna and Rechnitzer~\cite{MishnaRechnitzer2009} and Melczer and Mishna~\cite{MelczerMishna2014} using the iterated kernel method, a variant of the kernel method.

As we know, $\mathcal{S}=[[-1,1],[1,-1],[1,1]]$ is one of the singular models and its generating function is not D-finite. Now we rewrite the complete generating function of $\mathcal{S}$ into the following form:
\[
O(x,y;t)=\sum_{\substack{n_1,n_2,n_3\geq 0 }}o(n_1,n_2,n_3)x^{ -n_1+n_2+n_3}y^{n_1-n_2+n_3}t^{n_1+n_2+n_3},
\]
where $o(n_1,n_2,n_3)$ denotes the number of walks in the quarter plane with the $i$-th element of $\mathcal{S}$ appears $n_i$ times,  $(-n_1+n_2+n_3,\ n_1-n_2+n_3)$ denotes the ending point. Suppose $\mathcal{S^{'}}$ is a multi-set which can be reduced to $\mathcal{S}$ through getting rid of the repeated steps, and the $i$-th element of $\mathcal{S}$ repeats $r_i$ times in $\mathcal{S^{'}}$. Then the generating function for $\mathcal{S^{'}}$ can be given as
\[
O^{'}(x,y;t)=\sum_{\substack{n_1,n_2,n_3\geq 0}}r_1^{n_1}r_2^{n_2}r_3^{n_3}o(n_1,n_2,n_3)x^{ -n_1+n_2+n_3}y^{n_1-n_2+n_3}t^{n_1+n_2+n_3}.
\]
It's easy to verify that
\[
O^{'}(x,y;t)=O(\sqrt{\frac{r_3}{r_1}}x,\sqrt{\frac{r_3}{r_2}}y;\sqrt{r_1r_2}t),
\]
which implies that $O^{'}(x,y;t)$ is not D-finite, since algebraic substitution does not change the D-finite property. There are $7$ of the $43$ singular models can be reduced to $\mathcal{S}=[[-1,1],[1,-1],[1,1]]$ or it's $x/y$ symmetry, and the above discussions show the corresponding generating function for these $7$ models are all not D-finite.

By similar discussions for another singular model $[[-1,1],[1,-1],[0,1]]$, one can prove the generating functions for another $11$ models are all non-D-finite.

Thus, we have shown that the generating functions of the excursions of the $366$ nonsingular models are all non-$D$-finite and $18$ singular models are with non-D-finite generating functions. According to this fact and results of~\cite{KurkovaRaschel2011,BostanRaschelSalvy2014,MelczerMishna2014,MishnaRechnitzer2009}, we conjecture that the generating functions of the left $43$ singular models are all non-$D$-finite.

\vspace{.2cm} \noindent{\bf Acknowledgments.}
 We wish to thank Professor Manuel Kauers for helpful suggestions. This work was supported by the 973 Project, the PCSIRT Project of the Ministry of Education and the National Science Foundation of China.

{\section*{Appendix}\label{SECT:appendix}

\begin{longtable}{|c|l|c|c}
\hline \rule{0pt}{5pt}
Numbers & Models & Reduced Models \\[3pt]
\hline \rule{0pt}{5pt}
$1$
&
$[[-1, 1], [1, -1], [1, 1]]$
&
\\[3pt]
$2$
&
$[[-1, 1], [-1, 1], [1, -1], [1, 1]]$
&
\\[3pt]
$3$
&
$[[-1, 1], [1, -1], [1, 1], [1, 1]]$
&
\\[3pt]
$4$
&
$[[-1, 1], [-1, 1], [-1, 1], [1, -1], [1, 1]]$
&
$[[-1, 1], [1, -1], [1, 1]]$
\\[3pt]
$5$
&
$[[-1, 1], [-1, 1], [1, -1], [1, -1], [1, 1]]$
&
\\[3pt]
$6$
&
$[[-1, 1], [-1, 1], [1, -1], [1, 1], [1, 1]]$
&
\\[3pt]
$7$
&
$[[-1, 1], [-1, 1], [1, -1], [1, 1], [1, -1], [1, 1]]$
&
\\[3pt]
\hline
$8$
&
$[[-1, 1], [0, 1], [1, -1]]$
&
\\[3pt]
$9$
&
$[[-1, 1], [-1, 1], [0, 1], [1, -1]]$
&
\\[3pt]
$10$
&
$[[-1, 1], [-1, 1], [1, -1], [1, 0]]$
&
\\[3pt]
$11$
&
$[[-1, 1], [0, 1], [0, 1], [1, -1]]$
&
\\[3pt]
$12$
&
$[[-1, 1], [-1, 1], [-1, 1], [0, 1], [1, -1]]$
&
\\[3pt]
$13$
&
$[[-1, 1], [-1, 1], [-1, 1], [1, -1], [1, 0]]$
&
$[[-1, 1], [1, -1], [0, 1]]$
\\[3pt]
$14$
&
$[[-1, 1], [-1, 1], [0, 1], [0, 1], [1, -1]]$
&
\\[3pt]
$15$
&
$[[-1, 1], [-1, 1], [0, 1], [1, -1], [1, -1]]$
&
\\[3pt]
$16$
&
$[[-1, 1], [-1, 1], [1, -1], [1, 0], [1, 0]]$
&
\\[3pt]
$17$
&
$[[-1, 1], [-1, 1], [-1, 1], [0, 1], [0, 1], [1, -1]]$
&
\\[3pt]
$18$
&
$[[-1, 1], [-1, 1], [0, 1], [0, 1], [1, -1], [1, -1]]$
&
\\[3pt]
\hline
$19$
&
$[[-1, 1], [0, 1], [1, -1], [1, 0]]$
&
\\[3pt]
$20$
&
$[[-1, 1], [-1, 1], [0, 1], [1, -1], [1, 0]]$
&
\\[3pt]
$21$
&
$[[-1, 1], [0, 1], [0, 1], [1, -1], [1, 0]]$
&
\\[3pt]
$22$
&
$[[-1, 1], [-1, 1], [-1, 1], [0, 1], [1, -1], [1, 0]]$
&
\\[3pt]
$23$
&
$[[-1, 1], [-1, 1], [0, 1], [0, 1], [1, -1], [1, 0]]$
&
$[[-1, 1], [1, -1], [0, 1], [1, 0]]$
\\[3pt]
$24$
&
$[[-1, 1], [-1, 1], [0, 1], [1, -1], [1, -1], [1, 0]]$
&
\\[3pt]
$25$
&
$[[-1, 1], [-1, 1], [0, 1], [1, -1], [1, 0], [1, 0]]$
&
\\[3pt]
$26$
&
$[[-1, 1], [0, 1], [0, 1], [1, -1], [1, 0], [1, 0]]$
&
\\[3pt]
\hline
$27$
&
$[[-1, 1], [0, 1], [1, -1], [1, 1]]$
&
\\[3pt]
$28$
&
$[[-1, 1], [-1, 1], [0, 1], [1, -1], [1, 1]]$
&
\\[3pt]
$29$
&
$[[-1, 1], [-1, 1], [1, -1], [1, 0], [1, 1]]$
&
\\[3pt]
$30$
&
$[[-1, 1], [0, 1], [0, 1], [1, -1], [1, 1]]$
&
\\[3pt]
$31$
&
$[[-1, 1], [0, 1], [1, -1], [1, 1], [1, 1]]$
&
\\[3pt]
$32$
&
$[[-1, 1], [-1, 1], [-1, 1], [0, 1], [1, -1], [1, 1]]$
&
\\[3pt]
$33$
&
$[[-1, 1], [-1, 1], [-1, 1], [1, -1], [1, 0], [1, 1]]$
&
$[[-1, 1], [1, -1], [0, 1], [1, 1]]$
\\[3pt]
$34$
&
$[[-1, 1], [-1, 1], [0, 1], [0, 1], [1, -1], [1, 1]]$
&
\\[3pt]
$35$
&
$[[-1, 1], [-1, 1], [0, 1], [1, -1], [1, -1], [1, 1]]$
&
\\[3pt]
$36$
&
$[[-1, 1], [-1, 1], [0, 1], [1, -1], [1, 1], [1, 1]]$
&
\\[3pt]
$37$
&
$[[-1, 1], [-1, 1], [1, -1], [1, 0], [1, 0], [1, 1]]$
&
\\[3pt]
$38$
&
$[[-1, 1], [-1, 1], [1, -1], [1, 0], [1, 1], [1, 1]]$
&
\\[3pt]
$39$
&
$[[-1, 1], [0, 1], [0, 1], [1, -1], [1, 1], [1, 1]]$
&
\\[3pt]
\hline
$40$
&
$[[-1, 1], [0, 1], [1, -1], [1, 0], [1, 1]]$
&
\\[3pt]
$41$
&
$[[-1, 1], [-1, 1], [0, 1], [1, -1], [1, 0], [1, 1]]$
&
$[[-1, 1], [1, -1], [0, 1], [1, 0], [1, 1]]$
\\[3pt]
$42$
&
$[[-1, 1], [0, 1], [0, 1], [1, -1], [1, 0], [1, 1]]$
&
\\[3pt]
$43$
&
$[[-1, 1], [0, 1], [1, -1], [1, 0], [1, 1], [1, 1]]$
&
\\[3pt]
\hline
\caption{$43$ singular models.}\label{TB:43}

\end{longtable}

\end{document}